\documentclass[11pt]{amsart}
\usepackage[utf8]{inputenc}   
\usepackage[T1]{fontenc}
\usepackage{amsmath}     
\usepackage{amssymb}
\usepackage{mathrsfs}
\usepackage{dsfont}        
\usepackage{bbm}
\usepackage{amsthm}
\usepackage[ngerman, english]{babel}            
\usepackage{graphicx}
\usepackage[english] {babel}
\usepackage[a4paper,margin=3cm]{geometry}
\usepackage{graphicx} 
\usepackage{tikz}

\usepackage[all]{xy}  

\usepackage{color}

\usepackage[toc,page,title]{appendix}  
\usepackage{hyperref}
\hypersetup{
    colorlinks=true,
    linkcolor=blue,
    filecolor=magenta,      
    urlcolor=cyan,
  }

\renewcommand\Im{\operatorname{Im}}
\renewcommand\Re{\operatorname{Re}}

\DeclareMathOperator{\tr}{tr}

\DeclareMathOperator{\R}{\mathbb{R}}
\DeclareMathOperator{\N}{\mathbb{N}}

\DeclareMathOperator{\C}{\mathbb{C}}
\DeclareMathOperator{\Z}{\mathbb{Z}}

\DeclareMathOperator{\HH}{\mathcal{H}}

\DeclareMathOperator{\WW}{\mathcal{W}}
\DeclareMathOperator{\VV}{\mathcal{V}}

\DeclareMathOperator{\FF}{\mathcal{F}}
\DeclareMathOperator{\AAA}{\mathcal{A}}
\DeclareMathOperator{\TT}{\mathcal{T}}

\DeclareMathOperator{\A}{\mathring{A}} 

\DeclareMathOperator{\Ric}{Ric}
\DeclareMathOperator{\Sc}{Sc}
\DeclareMathOperator{\RR}{R}

\DeclareMathOperator{\diam}{diam}
\renewcommand{\d}{\mathrm{d}}
\DeclareMathOperator{\p}{\partial}
\DeclareMathOperator{\Hess}{Hess}
\DeclareMathOperator{\Vol}{Vol}

\newcommand{\suchthat}{\,\middle| \, }

\theoremstyle{plain}
\newtheorem{theorem}{Theorem}[section]
\newtheorem{lemma}[theorem]{Lemma}
\newtheorem{prop}[theorem]{Proposition}
\newtheorem{corollary}[theorem]{Corollary}

\theoremstyle{definition}
\newtheorem{definition}[theorem]{Definition}

\theoremstyle{remark}
\newtheorem{remark}[theorem]{Remark}
\newtheorem*{note}{Note}

\numberwithin{equation}{section}  

\title{Minimizers of Generalized Willmore Functionals}
\author{Alexander Friedrich } 
\address{University of Potsdam, Institute of Mathematics, Karl-Liebknecht-Str. 24-25,
D-14476 Potsdam OT Golm}
\thanks{The author was supported by the DFG project ME3816/1-2. Further, the author would like to thank Jan Metzger for his guidance and patience during the authors Ph.D. from which this article developed. Additional thanks goes to Christian Scholz and Michael Schwarz for a short but fruitful discussion on graph coloring.}

\date{September 5, 2019}

\begin{document}

\begin{abstract}
We introduce a notion of generalized Willmore functionals motivated by the Hawking energy of General Relativity and bending energies of membranes.
An example of a bending energy is discussed in detail.
Using results of Y. Chen and J. Li, we present a compactness result for branched, immersed, haunted, stratified surface with bounded area and Willmore energy.
This allows us to prove the existence of area constrained minimizers for generalized Willmore functionals in the class of haunted, branched, immersed  bubble trees by direct minimization.
Here a haunted, stratified surfaces are introduced, in order to account for bubbling and vanishing components along the minimization process. 
Similarly, we obtain the existence of area and volume constrained, minimal, closed membranes for the discussed bending energy.
Moreover, we argue that the regularity results of A. Mondino and T. Rivi\`ere for Willmore surfaces can be carried over to the setting of generalized Willmore surfaces. 
In particular, this means that critical points of a generalized Willmore functional are smooth away from finitely many points.
\end{abstract}

\thispagestyle{empty}
 
\maketitle

\pagestyle{headings}

\section{Introduction}
In this paper we study generalizations to the Willmore functional. 
That is we propose a notion of generalized Willmore functional (see Definition \ref{def_gen_will}) which is inspired by the Hawking energy of General Relativity and bending energies for thin membranes.

Let $\Sigma$ be an immersed, oriented surfaces in an oriented n-dimensional Riemannian manifold $(M^n,g)$.
 The metric $g$ induces a metric $\gamma$ on $\Sigma$. 
The second fundamental form of $\Sigma$ in $M$ defined by 
\[\vec{A}(X,Y):= -(\nabla^M_X Y)^{\perp}, \]
and the mean curvature vector $\vec{H}:= \tr_\Sigma \vec{A}$ is the  trace of $\vec{A}$. It allows us to define the Willmore energy:
\begin{align*}
\WW[\Sigma] &:= \frac{1}{4} \int_\Sigma |\vec{H}|^2 d\mu .
\end{align*}
In the example and the regularity discussion we will restrict to one codimension in which case we choose an outward pointing unit normal vector field $\nu$ and decompose $\vec{A} = A \otimes \nu  $ as well as $\vec{H} = H \nu$. We denote the area of $\Sigma$ by $|\Sigma|$.

 The Hawking energy is a quasi local energy given by
\[ \mathcal{E}[\Sigma] = \sqrt{\frac{|\Sigma|}{16 \pi}} \left( 1 - \frac{1}{16 \pi} \int_\Sigma H^2 - P(x,\nu)^2 \, \d \mu \right). \]
Here $\Sigma$ is a spherical surface in a three dimensional ambient Riemannian manifold $M$.
In General Relativity $M$ is embedded in a four dimensional Lorentz manifold $L$ with second fundamental form $K$ and $P=\tr_\Sigma K$. 
Clearly, we can analyze $\mathcal{E}$ subject to area constraint by investigating $\HH[\Sigma] = \int_\Sigma H^2 - P^2 \, \d \mu$ subject to area constraint. 
In \cite{Me_concentration_small_surf} we study this type of functionals in greater detail. Using the methods from \cite{LMI} and \cite{LMII}, we identify points in the ambient manifold around which small spherical minimizers concentrate. These can be interpreted as concentration points of the energy density as seen by the Hawking energy. Moreover, we calculate its expansion on coordinate spheres, which again reflects the concentration points.

 The Helfrich model represents thin membranes as surfaces which are a critical points of a bending energy under area and volume constraints. 
For a constant $c$, the following energy is proposed in \cite{Helfrich_elastic_properties}.
\[ \HH[\Sigma] = \int_\Sigma (H+c)^2 \, \d\mu \] 
We show that the main ideas of the theory of Willmore surfaces remain applicable for generalized Willmore functionals. 
In particular, we show the existence of minimizers through direct minimization and prove regularity away from finitely many points.
In the context of membranes this is novel since the existing literature usually  poses additional symmetry assumption for the surface which reduce the Euler-Lagrange equation to an ODE. 
Additionally, we have the advantage to work in the context of stratified surfaces directly. 
Hence, budding of membranes poses no problems for the framework.

In \cite{Chen_Li_14}  J. Chen and Y. Li investigated stratified surfaces and bubbling in the context of Willmore surfaces. 
We adopt their ideas and prove a compactness theorem for stratified surfaces.
In the following, a bubble forest is a stratified surface consisting of a closed Riemann surface with finitely many bubble trees attached. 
Haunted immersions are introduced in order to deal with degenerating tree structure of the bubble forest. They allow for some bubbles to be mapped to a point.

Our main result is  roughly the following, heuristic theorem. For the proper statement see Theorem \ref{thm_bubble_forest_cpt} along with the definitions of Section \ref{sec_preliminaries} and Definition \ref{def_haunted_immersion}.
\begin{theorem}
Let $S^k$ be a sequence of compact bubble forests. Let $\phi_k \in W^{2,2}(S^k, \R^n)$ be a sequence of irreducible, haunted, branched conformal immersions.
Assume $\phi_k$, the area $\AAA[\phi_k] $  and $ \WW[\phi_k]$ are uniformly bounded. 

Then $ \phi_k(S^k)$ either converges to a point or subconverges to an immersed haunted bubble tree $\phi(S)$. In the second case we find
\begin{align*}
\AAA[\phi] &= \lim_{k\to \infty} \AAA[\phi_k], \\
\WW[\phi] & \leq \lim_{k\to \infty} \WW[\phi_k] .
\end{align*}
\end{theorem}

 In order to employ this theorem in the search for minimizers we need to ensure that a given sequence does not shrink to a point. 
A straight forward way to establish this is to fix the area or a similar geometric quantity and solve the variational problem under constraints.
 Both of our examples are formulated as constrained variational problems and in fact our definition of generalized Willmore functionals is made such that this approach succeeds.
Hence we find the following existence results via direct minimization.
Here $\FF_a(\TT, M)$ denotes the space of haunted, branched immersions of bubble trees with area $a$ and the notion of $a$-generalized Willmore functional is introduced in Definition \ref{def_gen_will}.

\begin{theorem}
Let $(M,g)$ be compact Riemannian manifold and  let $\HH$ be an  $a$-generalized Willmore functional, then $ \inf \left\{ \HH[\phi] \mid \phi \in \FF_a(\TT, M) \right\}$ is attained in $\FF_a(\TT, M)$. 
\end{theorem}

\begin{corollary}
Let $(M,g)$ be a non compact Riemannian manifold with $C_B$ bounded geometry and  let $\HH$ be an $a$-generalized Willmore functional. Suppose there exists a transitive group action on $M$ that leaves $\HH$ and $\AAA$ invariant, then $ \inf \left\{ \HH[\phi] \mid \phi \in \FF_a(\TT, M) \right\}$ is attained in $\FF_a(\TT, M)$.
\end{corollary}

Similarly, we show the existence of minimal membranes with prescribed area and enclosed volume. 
The volume functional on $\FF(S, \R^3)$ reeds
\[ \VV[\phi] := \frac{1}{3}\int_{\Sigma} \langle x-x_0 , \nu \rangle \, \d \mu \]
and the bending engery in question is given by
\[\HH_{c, b}[\Sigma] := \int_\Sigma  (H+c)^2 \, \d \mu + b \left( \int_\Sigma H  \, \d \mu \right)^2,   \]
where $c$ and $b$ are constants. For $a,v \in \R^+$ define \[\FF_{a,v}(\TT, \R^3) :=\left\{  \phi \in \FF(\TT, \R^3) \mid \AAA[\phi] = a, \, \VV[\phi] = v \right\}. \] 
\begin{theorem}
For any $c, b \in \R$ and $a, v \in \R^+ $ such that $ 3 \sqrt{4 \pi} v \leq a^{3/2}$ and $-ab \leq 1$ the infimum of $\HH_{c,b}$ on  $\FF_{a,v}(\TT, \R^3)$  is attained.  
\end{theorem}

 Moreover, it is possible to use the topological and geometric structure of the target manifold to prevent shrinking. 
In \cite{Mon_Riv_13_1} A. Mondino and T. Rivi\`ere employ a curvature condition on the compact target to ensure the area of the of a spherical minimizer of $\AAA + \|\vec{A}\|^2_{L^2}$ stays bounded away from zero. 
In the same article they also obtain an area constrained Willmore minimizing bubble tree in every non trivial 2-homotopy class of the target. 
In \cite{Chen_Li_14} J. Chen and Y. Li construct minimizers of $\AAA + \WW$ in $S^n$, containing prescribed points. 
Additionally, they investigate Douglas type conditions which serve to exclude bubbling.

Finally, we turn to the regularity theory as developed by A. Mondino and T. Rivière in \cite{Mon_Riv_13_1} and show the following theorem in an analogous manner.
\begin{theorem}
Let $\phi\in W^{2,2} (D, M^3)$ be a conformal immersion with conformal factor $e^{2\lambda}$, $\lambda \in L^\infty(D)$. If $\phi$ solves $\Delta H + H|\A|^2 + F(\phi, d \phi, \nabla d \phi) =0$, where $F$ is as in Definition \ref{def_gen_will}, then $\phi$ is smooth.
\end{theorem}
In particular, this means that critical points of generalized Willmore functionals in codimension one are smooth, away from finitely many points. 
Although we prove the regularity only in the codimension one case we conjecture that it holds in any codimension. This is because the proof of the regularity of the Willmore equation in \cite{Mon_Riv_13_1} is carried out in arbitrary codimension.

The overview of the article is as follows. 
In Section \ref{sec_preliminaries} we introduce the definitions for branched, conformal immersions, stratified surfaces, bubble trees and their convergence, as well as generalized Willmore functionals, where we took inspiration from \cite{Chen_Li_14}, \cite{Kuwert_Li_conf_immersions} and \cite{Mon_Riv_13_1}. 
Section \ref{section_membranes} discusses bending energies for membranes as an example for generalized Willmore functionals, thus demonstrating that the theory developed here has immediate application in mathematical physics.
We introduce haunted immersions and prove our main results in Section \ref{sec_existence},  relying on \cite{Chen_Li_14} to show our compactness result in Theorem \ref{thm_bubble_forest_cpt}. The existence of minimizers then follows by direct minimization under constraints.
The final Section \ref{sec_regularity} is dedicated to the regularity theory for generalized Willmore equations.  Here we follow \cite{Mon_Riv_13_1} closely.


\section{Preliminaries}
\label{sec_preliminaries}

\begin{definition}[cf. {\cite[Definition 1 and 2]{Chen_Li_14}}]

Let $(S,\eta)$ be a Riemann surface and let $(M^n,g)$ be an n-dimensional, orientated Riemannian manifold which we assume to be isometrically embedded in some $\R^N$. 

\begin{enumerate}
\item For $k\in  \Z$ and $p \in [1,\infty]$ we define the Sobolev spaces as follows:
\[W^{k,p}(S,M) := \left\{ \phi\in  W^{k,p}(S, \R^N) \suchthat \phi(S) \subset M \textup{ a.e.} \right\}. \] 

\item
An element $\phi \in W^{2,2}(S,M)$ is called \emph{conformal immersion}, if $\phi$ is an immersion almost everywhere and if there is a function $e^{2\lambda}: S \to \R$, called the \emph{conformal factor} of $\phi$, such that
\[ \phi^* g = e^{2\lambda} \eta . \]

	\item  
	We say $\phi: S \to M$ is a \emph{branched conformal immersion} with finitely many branch points $ B \subset S$, if $\phi \in W^{2,2}_{loc}(S \setminus B ,M)$ is a conformal immersion and if for all $p \in B$ there is an open neighborhood $U_p$ and a constant $C$ such that
	\begin{align*}
	\int_{U_p \setminus \{p\}} 1+ |\vec{A}|^2 \, \d \mu_g  \leq C .
	\end{align*}
	
	\item Set 
	\begin{align*}
	\FF(S,M) &:= \{ \phi \in W^{2,2} (S,M) \mid \phi \textup{ is branched, conformal, immersion}\\
	& \hspace{0.8cm} \textup{with  branch points } B; \phi \in  W^{1,\infty}_{loc}(S\setminus B, M) \}
	\end{align*}
	and for $a>0$ define $\FF_a(S,M):= \{ \phi \in \FF(S,M) \mid |\phi(S)|= a \}$

\end{enumerate}
\end{definition}
Note that the conformal factor of a $\phi \in \FF(S,M)$ satisfies $\lambda \in L^{\infty}_{loc}$ away from the branch points.

 For an immersion $\phi \in \FF(S,M) $, $\Sigma:= \phi(S)$, we use $\WW[\Sigma]$ and $\WW[\phi]$ interchangeably. 
Moreover, at times we write $\AAA[\Sigma]$ or $\AAA[\phi]$ for the area $|\Sigma|$ in order to emphasize its role as a functional.
 E. Kuwert and Y. Li showed that branched conformal immersions can be extended to $W^{2,2}$ maps. 

\begin{theorem}[see {\cite[Theorem 3.1]{Kuwert_Li_conf_immersions}} ]
\label{thm_lift_sing}
Let $D$ be the unit disc in $\R^2$ and let $\phi \in W^{2,2}_{loc}(D\setminus \{0\}, \R^n)$, $n\geq 3$, be a conformal immersion, $\phi^*g =e^{2 \lambda } \delta$. If $\phi$ satisfies 
\[ \int_{U_p \setminus \{p\}} 1+ |\vec{A}|^2 \, \d \mu_g   \leq \infty \] 
then $\phi \in W^{2,2}(D,\R^n)$ and in complex coordinates we have
\begin{align*}
 \lambda (z) = m \ln |z| + w(z), \\
- \Delta  \lambda  = -2m \pi \delta_0 + K_g e^{2 \lambda } .
\end{align*}
Here, $m\in \N$, $w \in C^{0} \cap W^{1,2}(D)$, $K_g$ is the Gauss curvature of $g$ and $\delta_0$ is the delta distribution at $0$.
Additionally, the multiplicity of the immersion at $p=\phi(0)$ is given by
\[ \theta^2(\phi, p) = \# \phi^{-1}(p) = m+1 .\]
\end{theorem}

 The well known phenomenon of bubbling of $W^{2,2}$ immersions necessitates the introduction of stratified surfaces. 
\begin{definition}[ cf. {\cite[Definition 3]{Chen_Li_14}} ]
A compact connected metric space $(S,d)$ is called a \emph{stratified surface with singular points $P$}, if $P \subset S$ is a finite set such that:
\begin{enumerate}
	\item the regular part, $S \setminus P$, is a smooth Riemann surface without boundary. It carries a smooth metric $\eta$, whose induced distance function agrees with $d$.
 	\item  Moreover, for each $p\in P$ there is a $\delta > 0$ such that $B_\delta(p) \cap P = \{p\}$ and $B_\delta (p) \setminus \{p\} = \bigcup_{i=1}^{m(p)} \Omega_i$. Here $1<m(p) < \infty$ and the $\Omega_i$ are topological discs with one point removed. Additionally, we assume that $\eta$ can be extended to a smooth metric on each $\Omega_i \cup \{ p \}$.
\end{enumerate}
\end{definition}

 For a stratified surface, the regular part naturally decomposes into finitely many punctured connected Riemann surfaces  $S \setminus P = \bigcup_i S^i$. By the second point of the previous definition, we can add finitely many points to each $S^i$ in order to obtain a Riemann surface $\overline{S^i}$. This allows us to interpret a stratified surface as a collection of touching Riemann surfaces. 
\begin{figure}[ht]
\centering
\includegraphics[scale=0.5]{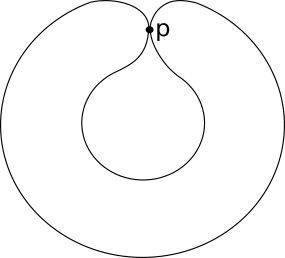}
\caption{A stratified torus with singular point $p$}
\label{fig_stratified_torus}
\end{figure}

 Consider the stratified torus $S$ in figure \ref{fig_stratified_torus}. It has one singular point $p$ and $S^1 = S\setminus \{p\}$ is a sphere with two punctures. We may add in two points $p_1, p_2$ such that $\overline{S^1}=S^1 \cup \{p_1\}\cup \{p_2\}$ is a sphere. In this picture we can understand $S$ as the immersion of $\overline{S^1}$.

By abuse of notation we usually denote a stratified surfaces as $S= \bigcup_i S^i$ and refer to Riemannian metrics on $S$ instead of on every $\overline{S}^i$.

\begin{definition} \
 \begin{enumerate}
\item 
Associate to every stratified surface $S = \bigcup_i S^i$ its \emph{dual graph}, where the vertices correspond to the components $S^i$ and two vertices are joined by an edge whenever the corresponding $S^i$  are joined by a singular point. 
Note that this construction allows for multiple edges and loops.

	\item 
A stratified surface whose regular part consists of punctured spheres and whose dual graph is a simple tree is called a \emph{bubble tree}. 
The constituting spheres are called \emph{bubbles}.

\item
If $S$ is a stratified surface and $S_1\subset S$ is a bubble tree, then we say $S_1$ is \emph{attached} to $S$ at $p\in S$ if $\overline{S\setminus S_1} \cap S_1 = \{ p \}$.

\item A stratified surface $S = S_0 \cup \bigcup_{i=1}^m S_i $ consisting of a Riemann surface $S_0$ with finitely many bubble trees attached at mutually distinct points is called \emph{bubble forest} with base $S_0$.
 Note that the dual graph of a bubble forest is still a tree.
\end{enumerate}
\end{definition}

\begin{figure}[ht]
\centering
\includegraphics[scale=0.5]{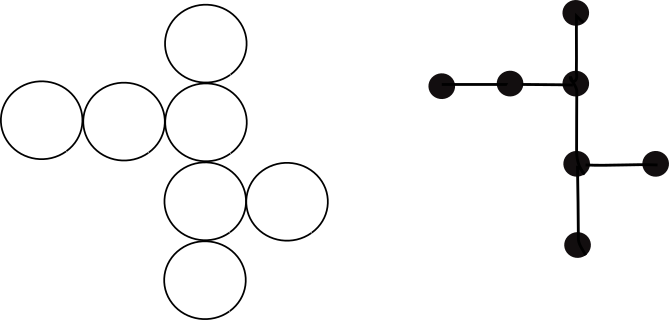}
\caption{A bubble tree and its dual graph}
\label{fig_bubble_tree}
\end{figure}

\begin{definition} 

\ \begin{enumerate}
	\item 
Let $S$ be a stratified surface with $S \setminus P = \bigcup_{i=1}^m S^i $ and let $M$ be a manifold of dimension three or higher. 
For $k\in \N$ and $p \in [1,\infty]$ denote by $W^{k,p}(S, M)$ the continuous maps $\phi: S \to M$ for which all $\phi|_{S^i}$ extend to  maps in $W^{k,p}(\overline{S^i},M)$.\\
 Additionally, we say that $\phi: S \to M$ is a (branched) immersion if all extensions $\phi|_{\overline{S^i}}$ are (branched) immersions.

\item 
Any functional defined for immersed surfaces we extend componentwise to immersed stratified surface. For example for $S$ and $\phi$ as above we set $\phi^i :=\phi|_{\overline{S^i}} $, $\Sigma = \phi(S)$ and $\Sigma^i := \phi^i(\overline{S^i})$; then the area and Willmore energy are given by
\begin{align*}
\AAA[\phi] &:= |\Sigma|:=  \sum_{i=1}^m |\Sigma^i| = \sum_{i=1}^m \AAA[\phi^i], \\
\WW[\phi] &:= \sum_{i=1}^m \WW[\phi^i] .
\end{align*}
\end{enumerate}
\end{definition}

\begin{definition}
\label{def_convergence_strat}
Let $(S,\eta_k)$ be a sequence of compact Riemann surfaces and $\phi_k\in$ \allowbreak $W^{2,2}(S,M)$ a sequence of branched, conformal immersions with conformal factors $e^{2 \lambda_k}$. Let $(S_\infty, \eta)$ be a stratified surface with singular set $P$ and let $\phi \in W^{2,2}(S_\infty,M)$ be a branched, conformal immersion with branch points $B$. We say $(S,\eta_k, \phi_k)$ converges to $(S_\infty, \eta, \phi)$ as immersed, stratified surfaces, if for all $k\in \N$ we can find open sets $U_k \subset S$ and $V_k \subset S_\infty$ such that

\begin{enumerate}

	\item $V_k \subset V_{k+1}$ and $P= S_\infty \setminus \bigcup_{k=1}^\infty V_k$. Moreover, $S_\infty \setminus V_k $ is a union of topological discs with finitely many smaller discs removed.
	
	\item $S\setminus U_k$ is a smooth surface with boundary and $\phi_k(S \setminus U_k)$ converges to $\phi(P)$ in Hausdorff distance.
	
	\item $\phi_k(S)$ converges to $\phi(S_\infty)$ in Hausdorff distance.
	
	\item There is a sequence of diffeomorphisms $\psi_k : V_k \to U_k$ such that $\phi_k \circ\psi_k \rightharpoonup \phi $ weakly in $W^{2,2}(K,M)$.
	
	\item The metrics $\psi_k^*\eta_k$ converge smoothly to  $\eta$.
\end{enumerate}
Further, let $(S= \bigcup_{i=0}^m S_i, \eta_k ) $ be a sequence of stratified surfaces  and $\phi_k \in W^{2,2}(S,M)$ a sequence of branched, conformal immersions.  
We say the sequence $(S,\eta_k, \phi_k)$ converges to $(S_\infty,\eta,\phi)$ as immersed, stratified surfaces if $(S_i, \eta_k|_{S_i}, \phi_k|_{S_i})$ converges to $(S_i^\infty, \eta_i, \phi^i)$ as immersed, stratified surfaces for all $i\in\{0,...,m\}$ and $S_\infty =\bigcup_{i=0}^m S_i^\infty $, $\eta|_{S_i^\infty} = \eta_i$ and $\phi|_{S_i^\infty} = \phi^i$. 
\end{definition}

\begin{figure}[ht]
\centering
\includegraphics[scale=0.3]{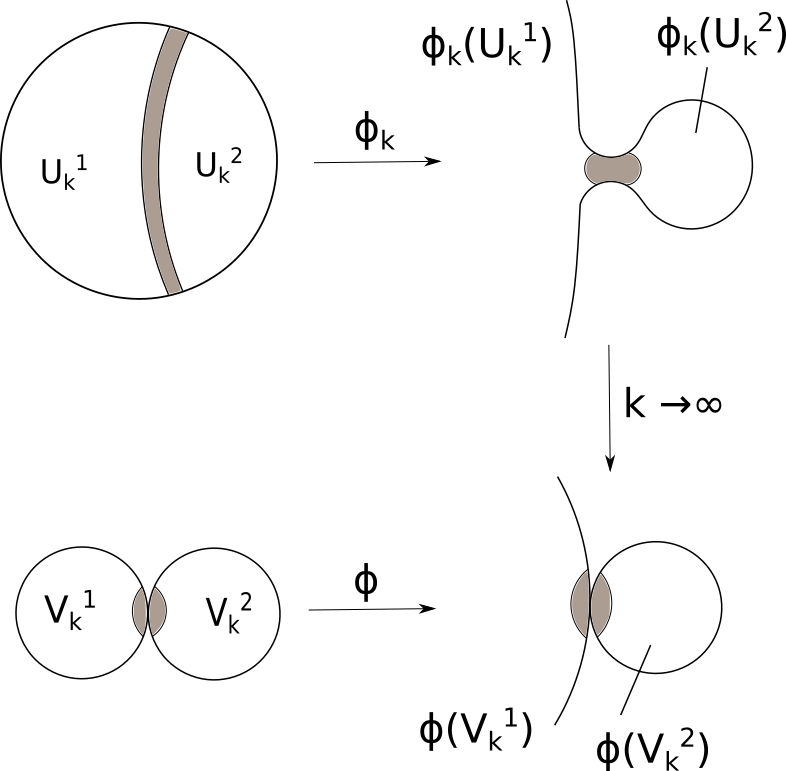} 
\caption{Side view of bubbling on a sphere}
\label{fig_bubbling}
\end{figure}

\begin{definition}
\label{def_gen_will}
Let $S$ be a closed stratified surface and let $(M,g)$  be an oriented n-dimensional Riemannian manifold. 
For $\phi \in \FF(S,M)$ denote the conformal factor by $e^{2\lambda}$ and the Hessian by  $\nabla d\phi$. 
Let $\{e_i\}_{i=1}^2$ be a local orthonormal frame on $TS$ and let $\{ \nu_i\}_{i=1}^{n-2}$ be local orthonormal frame of the normal bundle $N S$.
\ \begin{enumerate}
	\item A branched conformal immersion $\phi \in \FF(S, M)$ is said to solve a \emph{generalized Willmore equation} (away from the branch points) if
	\begin{align}
	\label{eq_gen_Will}
	\Delta_\perp \vec{H} + \sum_{i, j =1}^2 g(\vec{A}(e_i,e_j), \vec{H}) \vec{A}(e_i, e_j) - |\vec{H}|^2 \vec{H}+ F(\phi, d \phi, \nabla d\phi) = 0.
	\end{align}
	Here $\Delta_\perp$ denotes the Laplace operator on $N S$ and  $F(\phi, d \phi, \nabla d \phi): S \to N S$ is such that locally in conformal coordinates with $\lambda \in L^\infty$ and $F= F^i \nu_i$ we have
	\begin{align*}
	e^{2\lambda} F^i(\phi, d \phi, \nabla d\phi) &\in L^1 + W^{-1,2}(S) \ \ \textup{if } \phi \in \FF(S,M)  \\
	 e^{2\lambda}  F(\phi, d \phi, \nabla d\phi) &\in W^{k-1,l} (S, NS), \ l= \frac{2p}{2+p} + \epsilon \ \ \textup{if } \phi \in W^{k+2,p} , \ p>2 , \ k\geq 0
	\end{align*}
	for some $\epsilon >0$.

	\item A functional $ \HH$ on $\FF(S,M)$ is called an $a$-\emph{generalized Willmore functional}	if
\begin{enumerate}
	\item for any $\phi\in \FF_a(S,M)$ a bound $\HH[\phi] \leq \Lambda$ implies a bound on the Willmore energy $\WW[\phi] \leq C(\Lambda, a, M, \HH)$.
	
	\item $\HH$ is bounded from below on $\FF_a(S,M)$.
\item $\HH$ is invariant under diffeomorphisms of $S$.

	\item Let $\{\phi_k\}$ be a sequence in $\FF(S,M)$ with conformal factors $e^{2 \lambda_k}$. For any finite set $\mathfrak{S}\subset S$ the weak convergence $\phi_k \rightharpoonup  \phi$  in $W^{2,2}_{loc}(S\setminus \mathfrak{S},M)$ together with $ \|\lambda_k\|_{L^\infty(K)} \leq C_K$ for any $K \subset \subset S\setminus \mathfrak{S}$ implies $\HH[\phi] \leq \lim_{k\to \infty} \HH[\phi_k]$.

	\item $\HH$ is differentiable and its Euler-Lagrange equation is a generalized Willmore equation.
\end{enumerate}
If a functional is an $a$-generalized Willmore functional for all $a>0$ or if the area in question is understood we will simply refer to them as generalized Willmore functionals.
\end{enumerate}
\end{definition}

\begin{note} \ 
\begin{enumerate}

\item The generalized Willmore equation is of course based on the Euler-Lagrange equation of the Willmore functional which reads 
\[ \Delta_\perp \vec{H} + \sum_{i, j =1}^2 g(\vec{A}(e_i,e_j), \vec{H}) \vec{A}(e_i, e_j) - |\vec{H}|^2 \vec{H} -  \sum_{i=1}^2 \left(\RR^M(\vec{H},e_i)e_i \right)^\perp = 0 .\]

\item If a generalized Willmore equation is induced by a generalized Willmore functional, it will of necessity only depend on invariant quantities, that is $F(\phi, \nabla \phi, \nabla d\phi) = \tilde{F}(y, G, \vec{A})$, where $y\in \phi(S)$ and $G$ is the Gauss map of $\phi(S)$. If $\phi \in W^{k+2,p} \cap W^{1,\infty}$ is a branched, conformal immersion, we have,  away from the branch points, $\lambda \in L^\infty$, $e^{2\lambda}, G \in W^{k+1,p} \cap L^\infty  $ and $\vec{A} \in W^{k,p}$.
 
\item The area constrained variations of any generalized Willmore functional yields a generalized Willmore equation as well since for a normal variation we find $\delta_X \AAA[\Sigma] = \int_\Sigma g(\vec{H}, X) \,\d\mu $ and $\vec{H}$ obeys the conditions of the first part of Definition \ref{def_gen_will}.

\end{enumerate}

\end{note}

\noindent Recall also the uniformisation theorem for Riemann surfaces, we just need the following consequence.
\begin{theorem}[see {\cite[Chapter 1]{Hummel_Gromov_cptness}}]
\label{thm_uniformisation}
Let $(S, \eta)$ be a compact Riemann surface, then $S$ is conformal to a sphere, a torus or a surface of higher genus with constant Gauss curvature $1$, $0$ or $-1$ respectively. 
Moreover, if $S$ is a sphere then any two metrics are conformal and there is only one with Gauss curvature $1$. 
If $(S,\eta)$ is a torus it is conformal to $\C /(1, a+bi) $ where $-\frac{1}{2} < a\leq \frac{1}{2}$, $b \geq 0$, $a^2 + b^2 \geq 1$ and $a\geq 0$ provided $a^2+b^2 = 1$.
\end{theorem}


\section{Example: Thin Membranes} 
\label{section_membranes}

We model a membrane as a branched, immersed, stratified surface $\Sigma \subset (\R^3, \langle \cdot, \cdot \rangle)$ with a para\-metrization $\phi \in \FF(S, \R^3)$ that minimizes the following bending energy under area and volume constraints.
\[ \HH[\Sigma]  := \HH_{c, b}[\Sigma] := \int_\Sigma  (H+c)^2 \, \d \mu + b \left( \int_\Sigma H  \, \d \mu \right)^2   \]
Here, $c$, the spontaneous curvature, is a function on $\Sigma $, and  $b$ is a constant.
The first part of this energy corresponds to the one proposed by Helfrich and the second part is a non local generalization. See \cite{Ou-Yang_Tu_review} for a relatively recent review of membranes as elastic materials.

In the case that $\Sigma$ is the smooth boundary of a domain $\Omega$ we use the divergence formula to rewrite the volume of $\Omega$. Let $x$ be the position vector field in $\R^3$ and let $x_0 \in \R^3$, then $\Vol(\Omega) = \frac{1}{3} \int_\Sigma \langle x-x_0, \nu \rangle \,\d\mu $, independently of the choice of $x_0$. This motivates the introduction of the functional 
\[ \VV[\phi] := \frac{1}{3}\int_{\Sigma} \langle x-x_0 , \nu \rangle \, \d \mu \]
on $\FF(S,\R^3)$. It is still independent of the base point $x_0$, which also implies it is invariant under translations. This follows from the fact that any variation of $\tilde{\VV}[\phi] = \int_{\Sigma}  \langle x_0 , \nu \rangle \, \d \mu $ vanishes. 
This in turn relies on the divergence formula for vector fields on $\Sigma$ (see \cite[Section 2]{mantegazza1996}). Considering the variation of $\phi$  induced by scaling, $\phi_t(x) = (1+t) \phi(x)$, yields
\[ 0 = \left. \frac{d}{dt}\right|_{t=0} \tilde{\VV}[\phi_t] = \left. \frac{d}{dt}\right|_{t=0}  (1+t)^2 \tilde{\VV}[\phi] = 2 \tilde{\VV}[\phi] . \]
Hence, a volume constraint for $\phi \in \FF(S, \R^3)$ is defined to be a constraint on $\VV[\phi]$. 

\begin{prop}
\label{prop_membranes}
In the setting from above suppose that  $c : \R^3 \to \R $ is smooth and bounded. If $ - b a < 1$, then $\HH$ is an $a$-generalized Willmore functional. Its area and volume constrained Euler-Lagrange equation reads
\begin{align*}
\Delta H + H |\A|^2 + H\Ric(\nu,\nu)  &= \frac{1}{2} H c^2 + H^2 c + (H+ 2c) dc(\nu) - \tr_\Sigma \Hess^M c  -\frac{1}{2} \lambda H - \frac{1}{2} p\\
& \hspace{0,7cm} - \left( c+ b \int_\Sigma H \, \d \mu \right) (|\A|^2 + \Ric(\nu,\nu)) + \frac{1}{2} b H^2 \int_\Sigma H \, \d\mu.
\end{align*} 
Here $\lambda$ and $p$ are the Lagrange parameters for the area and the enclosed volume respectively.
\end{prop}

\begin{proof}
Let $S$ be a stratified surface with singular set $P$, let $\phi \in \FF(S,\R^3)$ and $\phi(S)=\Sigma$. Suppose that $\HH[\Sigma] \leq \Lambda$ and $|\Sigma| = a$, then
\begin{align*}
4 \WW[\Sigma] & \leq \Lambda - \int_\Sigma  2 H c + c^2 \, \d \mu -  b \left( \int_\Sigma H  \, \d \mu \right)^2 \\
& \leq \Lambda + C(c) a^{1/2}  \WW[\Sigma]^{1/2} - b \left( \int_\Sigma H  \, \d \mu \right)^2 .
\end{align*}
If $b\geq 0$, then we  omit the last term. Solving the quadratic inequality yields
\[ \WW[\Sigma] \leq C(\Lambda, c, a). \]
If $b < 0$, we estimate further:
\begin{align*}
4 \WW[\Sigma] & \leq \Lambda + C(c) a^{1/2}  \WW[\Sigma]^{1/2} + 4 |b| a  \WW[\Sigma] .
\end{align*}
If $|b| a < 1$ we absorb the last term to the left and find
\[ \WW[\Sigma] \leq C(\Lambda, c, b, a). \] 
Similarly, we obtain an estimate from blow. If $b\geq 0$ then $\HH$ is positive, so suppose $b < 0$ and $|b| a \geq 1 - \epsilon  $, for an $\epsilon \in (0,1)$.
\begin{align*}
\HH[\Sigma] &= \int_\Sigma  (H+c)^2 \, \d \mu + b \left( \int_\Sigma H  \, \d \mu \right)^2 \\
&\geq \int_\Sigma  2 H c + c^2 \, \d \mu + 4 (1- |b|a) \WW[\Sigma] \\
&\geq \left(1- \frac{1}{\epsilon} \right) \int_\Sigma  c^2 \, \d \mu  + 4 (1 -\epsilon - |b|a) \WW[\Sigma]\\
& \geq \left(1- \frac{1}{\epsilon} \right) C(c) a
\end{align*}

\noindent Since $H$ is a geometric quantity and $c$ is a function on $\R^3$, $\HH$ is invariant under re\-parametrizations.

To show that $\HH$ is lower semi continuous under weak $W^{2,2}_{loc}(S\setminus P,\R^3)$ convergence, let $\phi_k \in \FF(S, \R^3 )$ be a sequence of conformal maps with conformal factor $e^{2\lambda_k}$ and let $K \subset \subset S\setminus P$ such that 
\begin{align}
\phi_k & \rightharpoonup \phi \ \ \textup{weakly in } W^{2,2}(K, \R^3)  \nonumber \ \ \textup{and} \\
\label{eq_lambda_bound}
\| \lambda_k\|_{L^\infty(K)} & \leq C_K 
\end{align}
 We already know that the Willmore energy is lower semi continuous under these conditions, see \cite[Lemma A.8]{Mon_Riv_13_1}.
 Since the weak $W^{2,2}$ convergence implies strong $W^{1,p}$ convergence, we have that $\phi_k \to \phi $ and $\nabla \phi_k \to \nabla \phi$  pointwise almost everywhere. 
Dominated convergence and the smoothness of $c$ yields
\[ \int_K  c^2 \circ \phi_k \, \d \mu_{\phi_k}   \to \int_K  c^2 \circ \phi  \, \d \mu_{\phi}  .\] 
Weak $W^{2,2}$ convergence and the uniform upper and lower bound on the conformal factor \eqref{eq_lambda_bound} implies that $H_k e^{2\lambda_k} \to H e^{2\lambda}$ and $H_k e^{2\lambda_k} c\circ \phi_k \to H e^{2\lambda} c\circ \phi $ weakly in $L^2$. Testing with $1$ yields claim.

Let $\phi: I \times S \to \R^3$ be a normal variation of $\Sigma$ with $\left. \frac{\partial}{\partial s} \right|_0 \phi = f \nu$. The behavior of the geometric quantities is widely known, see for instance 
\cite[Section 7]{Huisken_Polden_Geo_Evo}.  
To derive the area and volume constrained Euler-Lagrange equation we need to calculate $\delta_f \HH = \lambda \delta_f \AAA + p \delta_f \VV$. \\
In particular, we find $\delta_f \VV[\phi] =\int_{\Sigma} f \, \d \mu$. In terms of generalized Willmore equations we have 
\begin{align*}
\Delta H + H |\A|^2  + F = 0
\end{align*}
for
\begin{align*}
 F &= H\Ric(\nu,\nu)  -\frac{1}{2} H c^2 - H^2 c - (H+2c) dc(\nu) + \tr_\Sigma \Hess^M c +\frac{1}{2} \lambda H + \frac{1}{2} p\\
& \hspace{0,7cm} + \left( c+ b \int_\Sigma H \, \d \mu \right) (|\A|^2 + \Ric(\nu,\nu)) - \frac{1}{2} b H^2 \int_\Sigma H \, \d\mu
\end{align*}
If $\phi \in W^{2,2}\cap W^{1,\infty}$ is a conformal parametrization of $\Sigma$, then $ |\nabla \phi|^2  F \in L^1$.

In terms of higher regularity, its worst term is of the form  $ |\nabla \phi|^2 c\circ \phi\, |\A|^2 \nu $. If $\phi \in W^{k+2,p}\cap W^{1,\infty}$, $k\geq 0$, $p> 2$, then, due to the Sobolev embeddings $W^{k+2,p} \hookrightarrow W^{k+1,q} \hookrightarrow C^{k,\alpha}$, for all $1\leq q< \infty$ and some $\alpha \in (0,1)$, we have $|\nabla \phi|^2 \nu \in W^{k+1,p} $ and $c \circ \phi \in W^{k+2,p} $. Due to the $|\A|^2$ part we find $|\nabla \phi|^2F \nu \in W^{k, p/2}$.
\end{proof}


\section{Existence of Generalized Willmore Surfaces}
\label{sec_existence}

In order to make use of the results on Willmore surfaces in Euclidean space, we  briefly recall how the Willmore energy of a closed surface $\Sigma \hookrightarrow M \hookrightarrow \R^N$ with respect to $\R^N$ is controlled by its Willmore energy with respect to $M$ and its area. Here we assume that the target  Riemannian manifold manifold $(M,g)$ has been isometrically embedded in some $\R^N$, $\dim M < N$. This can always be achieved via Nash embedding. \\
Introduce the second fundamental form and the mean curvature vector of $\Sigma$ in $\R^N$ as $\bar{A}$ and $\bar{H}$ respectively.  The second fundamental form of $M$ in $\R^N$ is denoted by $K$. We have $\bar{A}= \vec{A}  + K$ as well as $\bar{H} = \vec{H } + P $, for $P := \tr_\Sigma K$. Since $\Sigma$ is compact, we easily see
\begin{align}
\label{eq_Willmore_M_Rn}
\int_\Sigma |\bar{H}|^2 \, \d \mu \leq \int_\Sigma |\vec{H}|^2 \, \d \mu + \sup_{\Sigma} |P|^2  |\Sigma|.
\end{align}
 In \cite{Chen_Li_14} J. Chen and Y. Li proved a Gauss-Bonnet formula for closed branched conformal immersions. 
\begin{lemma}[see {\cite[Lemma 3.2]{Chen_Li_14}}]
\label{lemma_Gauss_Bonnet}
If $\phi \in \FF(S,\R^n)$ and $\phi (S) = \Sigma$ then
\[ \int_\Sigma \Sc_\Sigma \, \d \mu = 8 \pi (1-q(S)) + 4 \pi b  . \]
Here $q(S)$ is the genus of $S$ and $b$ is number of branch points counted with multiplicity. 
\end{lemma}

Integrating over the Gauss equation yields
\begin{align}
\label{eq_branch_point_bound} 
2 \pi b \leq \frac{1}{4}  \int_\Sigma |\bar{H}|^2 \, \d \mu - 4 \pi (1-q(\Sigma)).
\end{align}
In the same paper J. Chen and Y. Li proved a powerful compactness result for $W^{2,2}$ branched, conformal immersions, which is the heart of our existence results.

\begin{theorem}[see {\cite[Theorem 1]{Chen_Li_14}}]
\label{thm_bubble_tree_conv}
Let $(S, \eta_k)$ be a sequence of closed Riemann surfaces with metrics as given by Theorem \ref{thm_uniformisation} and let $ \phi_k \in W^{2,2}((S,\eta_k), \R^n) $ be a sequence of branched conformal immersions, for some $n> 2$.
If $\phi_k(S) \cap B_{R_0} \neq \emptyset$ for some $R_0>0$, and if there are positive constants $a$ and $\Lambda$ such that
\begin{align*}
\AAA[\phi_k] &\leq a   \\
\WW[\phi_k] & \leq \Lambda  
\end{align*}
 for all $k\in \N$, then $\Sigma_k$ either  converges to a point or there is a stratified surface $(S_\infty, \eta)$ and a branched, conformal immersion $\phi \in W^{2,2}(S_\infty, \R^n) $ such that a subsequence of $(S, \eta_k, \phi_k)$ converges to $(S_\infty, \eta, \phi$) in the sense of immersed, stratified surfaces.\\
Moreover, 
\begin{align*}
\AAA[\phi] &= \lim_{k\to \infty} \AAA[\phi_k] ,\\
\WW[\phi] & \leq \lim_{k\to \infty} \WW[\phi_k] .
\end{align*}

\end{theorem}
\begin{remark}
\label{remark_convergence_detail}
From the proof of the theorem we learn more about the the convergence and the structure of $S_\infty$. 
It is obtained by attaching finitely many bubble trees to a stratified surface $T$. 
This base stratified surface $T$ in turn is formed as the limit of the $(S,\eta_k)$ as nodal surfaces with possibly some bubbling at the nodal points. Additionally, if the  $\eta_k$ smoothly converge to a smooth metric $\eta$ on $S$, then $T = (S,\eta)$ and $S_\infty$ is a bubble forest with base $S$. 
In this case, if $U_k \subset S $ are the open sets guaranteed by the convergence as immersed, stratified surfaces, we have $ \lim_{k\to 0}|\phi_k(S\setminus U_k)| = 0$.

Furthermore, let $P$ be the singular points of $S_\infty$. We may assume that the branch points of the sequence $\phi_k\circ \psi_k$ converge to a finite set $\tilde{B}$ and that there is a finite set $\mathfrak{S} \subset S_\infty \setminus P$ such that the  conformal factors $e^{2\tilde{\lambda}_k}$ corresponding to $\phi_k\circ \psi_k $ obey $\| \tilde{\lambda}_k \|_{L^\infty(K\cap V_k)} \leq C_K$ for all $K \subset \subset S \setminus (P \cup \mathfrak{S} \cup \tilde{B})$; cf. \cite[Proof of Theorem 1, Page 30]{Chen_Li_14}. 

\end{remark}

Since the convergence as immersed, stratified surfaces leaves the class of surfaces, we need to formulate a compactness theorem for stratified surfaces. Ultimately, our goal is the minimization of a generalized Willmore functional over a surface $S_0$, hence we restrict ourselves to bubble forests with base $S_0$. The idea is then to apply Theorem \ref{thm_bubble_tree_conv} to every part of the bubble forest. Unfortunately, this means that parts of the forest can vanish, even though the whole forest cannot due to constraints. This would destroy the tree structure of the dual graph and hence leave the class of immersed bubble forests. To remedy this we introduce ghost bubbles and haunted immersions.

\begin{definition} \ 
\label{def_haunted_immersion}
\begin{enumerate}
	\item  Let $S = \bigcup_{i=1}^m S_i$ be a stratified surface and let $\phi:S \to M$ be a continuous map into a manifold $M$. We say $\phi$ is a \emph{haunted} immersion, if it is constant on some, but not all, components of $S$ and an immersion on the rest. A component $S_i$ is called a \emph{ghost} if $\phi|_{S_i}$ is constant, otherwise it is called \emph{regular}.

\item Let $(S= \bigcup_{i=1}^m S_i,\eta_k)$ be a sequence of compact, stratified surfaces and \\ 
$\phi_k\in W^{2,2}(S,M)$ a sequence of haunted, branched, conformal immersions. 
Let $(S_\infty, \eta)$ be a stratified surface and let $\phi \in W^{2,2}(S_\infty,M)$ be a haunted, branched, conformal immersion. 
We say $(S,\eta_k, \phi_k)$ converge to $(S_\infty, \eta, \phi)$ as haunted, immersed, stratified surfaces if 
\begin{enumerate}
	\item $(S_i, \phi_k|_{S_i})$  converges to a point $x_i$ for some but not all $i$, setting $S_i^\infty = S_i$, $\phi^i = x_i$ and 
	
	\item the remaining $(S_i, \phi_k|_{S_i})$ converge to $(S_i^\infty, \eta_i, \phi^i)$ as immersed stratified surfaces such that
	
	\item $S_\infty =\bigcup_{i=0}^m S_i^\infty $, $\eta|_{S_i^\infty} = \eta_i$ and $\phi|_{S_i^\infty} = \phi^i$. 

\end{enumerate}

\end{enumerate}

\end{definition}

 Suppose $\phi$ is a haunted immersion of a bubble forest $S$. 
If a ghost is connected to only one other component then we delete it.
If a ghost is connected to two other components, say $S_1$ and $S_2$ with common points $p_1$ and $p_2$ respectively, we delete it as well and identify $p_1$ and $p_2$. 
Repeating this process until there are no ghosts left or until every ghost is connected to three or more components yields a bubble forest $S'$ (possibly with a different base then $S$) and a haunted immersion $ \phi' $ which is given by $\phi$ on every component of $S'$. A tuple $(S',\phi')$ obtained that way is called \emph{irreducible}.

\begin{figure}[ht]
\centering
\includegraphics[scale=0.5]{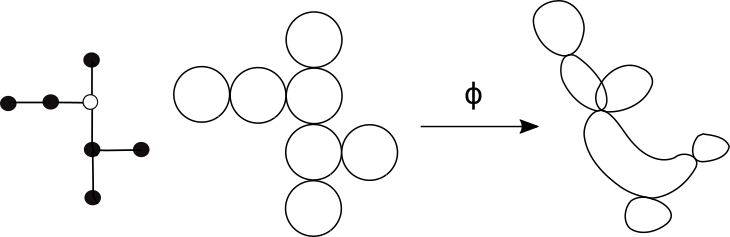}
\caption{An immersed haunted bubble tree and its dual graph, where the ghost is drawn white.}
\label{fig_haunted_bubble_tree}
\end{figure}
 The following lemma on graph coloring asserts that in an irreducible, haunted bubble forest the number of ghosts is bounded by the number of regular components.

\begin{lemma}
\label{lemma_graph_coloring}
Let $G$ be a finite tree, colored in black and white according to the rule: a vertex can be white if its degree is bigger or equal to three.  Then there are more black vertices then white ones. 
\end{lemma}

\begin{proof}
For a tree $G$ let $W(G)$ be the number of white vertices and $B(G)$ be the number of black vertices. Note that any endpoint of $G$ has to be black and the claim holds for trees with up to four vertices.

We argue by induction. Suppose the claim holds for trees with $n$ and $n-1$ vertices. Let $G$ be a tree with $n+1$ vertices and let $p$ be a boundary vertex connected to $q$.
\begin{enumerate}
	\item Suppose $\deg(q)=2$, then $q$ has to be black and $W(G)=W(G\setminus \{p\}) \leq B(G\setminus \{p\}) < B(G) $.
	\item Suppose $\deg(q) =3$ and $q$ is black or $\deg(q) \geq 4$, then there is no need to recolor $q \in G\setminus \{p\}$ and we have $W(G)=W(G\setminus \{p\}) \leq B(G\setminus \{p\}) < B(G) $ as before.
	\item Suppose $\deg(q) = 3$ and $q$ is white. Let $r$ and $s$ be the other vertices adjacent to $q$ and consider the tree $\tilde{G} = G\setminus \{p, q\}$ where we joined $r$ and $s$. Then
	\[ W(G)=  W(\tilde{G}) + 1 \leq B(\tilde{G}) + 1 = B(G). \]
\end{enumerate}
\end{proof}

\begin{theorem}
\label{thm_bubble_forest_cpt}
Let $(S^k,\eta_k)$ be a sequence of compact bubble forests with base $S_0$. Suppose $\eta_k$ is as given by Theorem \ref{thm_uniformisation} on every component of $S^k$, which fixes $\eta_k$ on every bubble. Suppose additionally that $\eta_k|_{S_0}$ converges smoothly to a smooth metric $\eta'$ on $S_0$.\\
Let $\phi_k \in W^{2,2}(S^k, \R^n)$ be a sequence of irreducible, haunted, branched, conformal immersions.
Assume  $\phi_k(S^k) \cap B_{R_0} \neq \emptyset$ for some $R_0$ and there are positive constants $a$ and $\Lambda$ such that 
\begin{align*}
\AAA[\phi_k] &\leq a \\
\WW[\phi_k] &\leq \Lambda 
\end{align*}
for all $k\in \N$.
Then there exists a bubble forest $S= S_0 \cup \bigcup_{i=1}^m S_i$,  a stratified surface $(\tilde{T} = \bigcup_{i=0}^{m'} T_i , \tilde{\eta}) $ and a haunted, branched, conformal immersion $\tilde{\phi} \in W^{2,2}(\tilde{T}, \R^n)$, 
a bubble forest $(T,\eta) $ (with base $S_0$ or a sphere) and an irreducible, haunted, branched, conformal immersion $\phi \in W^{2,2}(T, \R^n)$ such that $(\tilde{T}, \tilde{\eta}, \tilde{\phi})$  and $(T, \eta, \phi)$ differ only by ghosts and either 
\begin{enumerate}
	\item[1)] $\phi_k$ converges to a point, or
	
	\item[2)] there is, a subsequence of $\{\phi_k \}$ defined on $(S,\eta_k)$ such that $(S,\eta_k,\phi_k)$ converges to $(\tilde{T}, \tilde{\eta}, \tilde{\phi})$ as haunted immersed stratified surfaces.
\end{enumerate}
Moreover, if $U_k \subset S $ are the open sets guaranteed by the convergence as haunted, immersed, stratified surfaces, then $ \lim_{k\to 0}|\phi_k(S\setminus U_k)| = 0$ and
\begin{align*}
\AAA[\phi] &= \lim_{k\to \infty} \AAA[\phi_k], \\
\WW[\phi] & \leq \lim_{k\to \infty} \WW[\phi_k] .
\end{align*}
\end{theorem}

\begin{proof}
First, note that the number of  regular components of $S^k$ is uniformly bounded as they each have Willmore energy at least $4\pi$ and by Lemma \ref{lemma_graph_coloring} the total number of components is bounded. 
This means that there are only finitely many possible dual graphs along the sequence $S^k$ and we can choose subsequences of $S^k$ and $\phi_k$ such that they all agree. 
This means the $S^k$ agree as topological spaces but differ by their metric and their singular points $P^k$. 
Call the underlying topological space $S=S_0 \cup \bigcup_{i=1}^m S_i$ and note that the number of branch points is bounded by \eqref{eq_branch_point_bound}.

For an $i\in \{ 0,1,..., m \}$ consider $(S_i, \phi_k|_{S_i})$, we apply Theorem \ref{thm_bubble_tree_conv} so that either it becomes a ghost, setting $\tilde{S}_i = S_i$ and $\phi^i(x) = \lim_{k\to \infty} \phi_k(S_i)$,
 or $\left(S_i, \phi_k|_{S_i} \right)$ subconverges to $\left(\tilde{S}_{i}:= S_i \cup \bigcup_{j=1}^{m_i} S_{i,j}, \phi^i \right)$ as immersed stratified surfaces, where $\tilde{S}_i$ is a bubble forest and $\phi^i \in W^{2,2}(\tilde{S}_i, \R^n)$ is a branched, conformal immersion.

Next, we track the singular set of points $P^k$. For any $p \in P^k$ there are two components $S_i, S_j$ such that $p = S_i \cap S_j$, due to the tree structure of $S$. 

We may assume that $\phi_k(p)$ converges to a point $y \in \R^n$. Since $\phi^i(\tilde{S}_i)$ is compact, the distance $d(\phi_k(p), \phi(\tilde{S}_i))$ is attained for a sequence of points $y_k\in \phi^i(\tilde{S}_i) $. Now, $\phi_k(S_i)$ converges to  $\phi^i(\tilde{S}_i)$ in Hausdorff distance and we find $ y_k \to y \in \phi^i(\tilde{S}_i)$ as $\phi^i(\tilde{S}_i)$ is closed. Choose $x_i \in S_i$ such that $\phi^i(x_i)=y$.

Since the same reasoning holds for $\phi_k|_{S_j}$, we find a $x_j \in \tilde{S}_j$ such that $\phi^j(x_j) =y$. 
This means we can join the two bubble trees together: define $\tilde{T}_{ij} = \tilde{S}_i \sqcup \tilde{S}_j / (x_i \sim x_j)$ and $\phi^{ij} : \tilde{T}_{ij} \to \R^n$, $\phi^{ij}|_{\tilde{S}_i} := \phi^i$, $\phi^{ij}|_{\tilde{S}_j} := \phi^j$. 
In this way we join up all the bubble trees to obtain a haunted, branched, immersed, stratified surface $(\tilde{T} = \bigcup_{i=0}^m \tilde{S}_i, \tilde{\phi})$, where $\tilde{\phi}$ is given by $\phi^i$ on $\tilde{S}_i$ and is continuous throughout. 
It is then clear that $(S, \eta_k, \phi_k)$ converges to $(\tilde{T}, \tilde{\eta}, \tilde{\phi})$ as haunted, immersed, stratified surfaces. 
Here $\tilde{\eta} = \eta$ on every bubble and $\tilde{\eta} = \eta'$ on $S_0$, provided $\phi_k(S_0)$ does not shrink to a point. 

During this procedure it is possible to loose the tree structure, namely if two or more singular points of $\tilde{T}$ overlap which can only happen at the points constructed above. See figure \ref{fig_trublems} for illustration.

\begin{figure}[ht]
\centering
\includegraphics[scale=0.5]{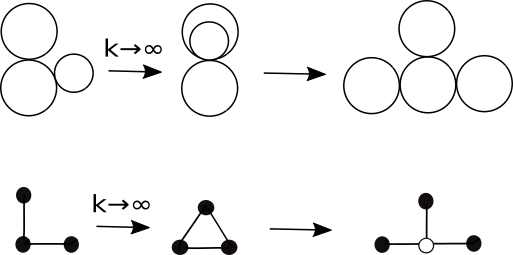}
\caption{Introducing ghosts into a degenerating bubble tree}
\label{fig_trublems}
\end{figure}

\noindent This is remedied by introducing a ghost. 
Let $p$ be a singular point of $\tilde{T}$ such that $p \in \bigcap_{j=1}^l \tilde{S}_{i_j}$ for $l>2$. Set $p_j := \{ p\} \cap \tilde{S}_{i_j}$, take a sphere $S'$ and $l$ mutually distinct points $\{a_j\}$ on it. 
Define the stratified surface $\tilde{T}' := \tilde{T} \sqcup S' /{\sim}$ where we no longer identify the $p_j$ but instead identify $p_j$ with $a_j$ and set $\tilde{\phi}'|_{\tilde{T}} = \tilde{\phi}$ and $\tilde{\phi}'|_{S'} = \tilde{\phi}(p)$. \\
Employing this method as often as necessary to obtain a tree and then deleting any unnecessary ghost yields the claim.

\end{proof}

 In terms of direct minimization the compactness result for haunted immersions of bubble forest puts the base in competition to the bubbles. 
Hence, we have to restrict to bubble trees, so as not to loose the base.

Let $\TT$ be the class of bubble trees. 
For a Riemannian manifold $(M,g)$ of dimension three or higher, a positive constant $a$ and an a-generalized Willmore functional $\HH$  set
\begin{align*} 
\FF(\TT,M)&:= \{ \phi \in \FF(S,M) \mid S \in \TT, \, (S,\phi) \textup{ is haunted} \} ,\\
\FF_a(\TT,M)&:= \{ \phi \in \FF(\TT,M) \mid \AAA[\phi]=a \} ,\\
\beta(\HH, M,a) &:= \inf \left\{ \HH[\phi]  \, \middle| \, \phi\in \FF_a(\TT,M) \right\} .
\end{align*}

\begin{theorem}
\label{theorem_existence_cpt}
Let $(M,g)$ be a compact Riemannian manifold and  let $\HH$ be an  $a$-generalized Willmore functional, then $\beta(\HH, M, a)$ is attained in $\FF_a(\TT, M)$. 
\end{theorem}

\begin{proof}
Pick a sequence $\phi_k \in \FF_a(\TT,M)$ realizing $ \beta(\HH, M, a)$. By definition, $\HH$ is bounded from below and along this sequence it is bounded from above. 
Again, by the definition of a-generalized Willmore functional, we know that the Willmore energy with respect to $M$ is bounded. This implies
\[ \WW[\phi_k, \R^N] \leq \WW[\phi_k, M] + C a \leq \Lambda (M, a, \beta, \HH ) .  \]
After reducing the $\phi_k$, if necessary, we are in the context of Theorem \ref{thm_bubble_forest_cpt}, as $M$ is compact. Since the area is fixed, the sequence $\phi_k$ cannot shrink to a point. 
The convergence as haunted, immersed, stratified surfaces  yields a limit $\phi \in \FF_a(\TT, M)$. 
By Remark \ref{remark_convergence_detail} and Definition \ref{def_gen_will} we know that $\HH$ is lower semi continuous with respect to this convergence, so we find 
\[ \beta(\HH, M, a) \leq \HH[\phi] \leq \lim_{k\to \infty} \HH[\phi_k] =  \beta(\HH, M, a).  \]
\end{proof}

\begin{corollary}
\label{cor_existence_homogeneous}
Let $(M,g)$ be a non compact Riemannian manifold with $C_B$ bounded geometry and  let $\HH$ be an a-generalized Willmore functional. Suppose there exists a transitive group action on $M$ that leaves $\HH$ and $\AAA$ invariant, then $\beta(\HH, M, a)$ is attained in $\FF_a(\TT, M)$.
\end{corollary}

\begin{proof}
Pick a sequence $\tilde{\phi}_k \in \FF_a(\TT,M)$ realizing $ \beta(\HH, M, a)$ and use the transitive action to obtain a sequence $\phi_k \in \FF_a(\TT,M)$ that still realizes $\beta(\HH, M, a)$ and whose images intersect a fixed point.
If $\Sigma_k$ is the image of $\phi_k$ then \cite[Lemma 2.5]{LMII} asserts
\[ \diam_M(\Sigma_k) \leq C(C_B) \left(|\Sigma_k|^{1/2} \WW[\Sigma_k, M]^{1/2} + |\Sigma_k| \right) \leq C\left(C_B, M, a, \beta, \HH \right) . \]
Since the diameter in $M \subset \R^N$ is larger then the one in $\R^N$, we conclude with Theorem \ref{theorem_existence_cpt}.
\end{proof}

\begin{theorem}
Let $\HH_{c,b}$ and $\VV$ be the functionals on $\FF(\TT, \R^3)$ introduced in Section \ref{section_membranes}; representing bending energy and enclosed volume. 
For $a,v \in \R^+$ define \[\FF_{a,v}(\TT, \R^3) :=\left\{  \phi \in \FF(\TT, \R^3) \mid \AAA[\phi] = a, \, \VV[\phi] = v \right\}. \] 
For any $c, b \in \R$ and $a, v \in \R^+ $ such that $ 3 \sqrt{4 \pi} v \leq a^{3/2}$ and $-ab \leq 1$ the infimum of $\HH_{c,b}$ on  $\FF_{a,v}(\TT, \R^3)$  is attained.  
\end{theorem}

\begin{proof}
Let $\{ \tilde{\phi}_k \}_{k\in \N}$ be a sequence in $\FF_{a,v}(\TT,  \R^3)$ realizing the infimum. 
Let $\{ T_k \}_{k\in \N}$ be a sequence of translations such that the image of $\phi_k := T_k \circ \tilde{\phi}_k$ contains the origin. 
Since $\HH_{c,b}$, $\AAA$ and $\VV$ are invariant under translations, $\HH_{c,b}[\phi_k]$ still converges to the infimum in $\FF_{a,v}(\TT, \R^3)$. 
If $-ab \leq 1$ then Proposition \ref{prop_membranes} asserts that $\HH_{c,b}$ is an $a$-generalized Willmore functional. 
As in the proof of Corollary \ref{cor_existence_homogeneous} we know that $\Im(\phi_k) \subset B_{R_0}(0)$
for a $R_0$ and all $k\in \N$. 
By Theorem \ref{thm_bubble_forest_cpt} we obtain a subsequence, again denoted by $\{ \phi_k \}$, $\phi_k \in \FF_{a,v}(S,\R^3)$, $S \in \TT$, that converges to $\phi \in \FF_a(T, \R^3)$ as haunted, immersed, stratified surfaces; where $(T, \phi)$ differs from a bubble tree only by ghosts. 
Moreover, we have
\[\HH_{c,b}[\phi] \leq \inf\{ \HH_{c,b}[\psi] \mid \psi \in \FF_{a,v}(\TT,\R^3) \}. \] 
Now we argue that $\VV[\phi] =v$.

Since the estimate $\VV\left[\phi_k|_{S^i}\right] \leq C \diam\left(\phi_k(S^i) \right) |\phi(S^i)|$ holds on any component $S^i$ of $S$, the bubbles that shrink to a point do not hold any volume in the limit.
Similarly, ghosts do not carry any volume, hence we will disregard them in the following.
Let $P$ be the set of singular points of $T$.
The convergence as haunted, immersed stratified surfaces yields, the existence of open sets $U_k \subset S$, $V_k \subset T$, $V_{k} \subset V_{k+1}$, $T\setminus \bigcup_{k\in \N} V_k \subset P$, and diffeomorphisms $\psi_k : V_k \to U_k$. 
Furthermore, we know $|\phi_k(S \setminus U_k)| \to 0$, $|\phi_k \circ \psi_k (V_k)| \to a$ and $|\phi(V_k)| \to a$. 
Let $\bigcup_{i=1}^m S^i$ be the union of all the components of $S$ such that $\phi_k(S^i)$ converges to a point. Set $\phi_k':=\phi_k \circ \psi_k$.
\begin{align*}
3 \left|v - \VV[\phi] \right| &\leq \left|  \int_{\phi_k'(V_k)} \langle x , \nu_k \rangle \,\d \mu_k - \int_{\phi(V_k)} \langle x , \nu \rangle \,\d \mu \right| + \sup_{x \in \phi_k(S)} |x| |\phi_k(S\setminus U_k)| \\
& \hspace{0.7cm}+ \sup_{x \in \phi(T)} |x| |\phi(T \setminus V_k)| + \left| \int_{\phi_k\left( \bigcup_{i=1}^m S^i \right)} \langle x , \nu_k \rangle \,\d\mu_k\right|
\end{align*}
For $j\in \N$ fixed we have 
\[ \int_{\phi_k'(V_j)} \langle x , \nu_k \rangle \,\d \mu_k \to \int_{\phi(V_j)} \langle x , \nu \rangle \,\d \mu \]
and
\[ \left|\phi_k'(V_j) \right| \to |\phi(V_j) |\]
 for $k \to \infty$ since $\phi_k'$ converges to $\phi$ weakly in $\WW^{2,2}_{loc}(T\setminus P, \R^3)$.
Additionally, we have 
\begin{align*}
\left| \int_{\phi(V_k)} \langle x , \nu \rangle \,\d \mu - \int_{\phi(V_j)} \langle x , \nu \rangle \,\d \mu \right| & \leq R_0|\phi(T) \setminus \phi(V_j)|,
\end{align*}
as well as
\begin{align*}
\left| \int_{\phi_k'(V_k)} \langle x , \nu_k \rangle \,\d \mu_k - \int_{\phi_k'(V_j)} \langle x , \nu_k \rangle \,\d \mu_k \right| & \leq  R_0 |\phi_k'(V_k) \setminus \phi_k'(V_j)|.
\end{align*}
For $\epsilon > 0 $ choose $j \in \N$ such that $|\phi(T) \setminus \phi( V_j)| \leq \frac{\epsilon}{4 R_0}$ and estimate
\begin{align*}
|\phi_k'(V_k) \setminus \phi_k'(V_j)| & \leq \left| | \phi_k'(V_k)| - a \right| + \left| a - |\phi(V_j)| \right| + \left| |\phi(V_j)| - |\phi_k'(V_j) | \right|.
\end{align*}
Now we may choose $k>j$ such that $|v- \VV[\phi]|\leq \epsilon$. 
\end{proof}


\section{Regularity of Generalized Willmore Surfaces}
\label{sec_regularity}

In this section we present the regularity theory for immersions of generalized Willmore type analogous to the theory developed by A. Mondino and T. Rivière in \cite{Mon_Riv_13_1} for critical points of the Willmore functional under conformal constraint. 
It hinges on the fact that the corresponding Euler-Lagrange equation exhibits a divergence form and that suitably chosen potentials obey the Laplace equation  with a Wente type structure. 

We cannot cite the regularity result directly as the function $F$ in Equation \eqref{eq_gen_Will} represents a more general nonlinearity then the ones treated in \cite{Mon_Riv_13_1}, though it turns out that we can follow the same arguments.

Consider a conformal embedding $\phi \in W^{2,2} \cap W^{1,\infty}(D,M)$ from the two dimensional open disc $(D, \langle \cdot , \cdot \rangle_E)$ to a three dimensional, oriented  Riemannian manifold $(M,\langle \cdot , \cdot \rangle)$ with conformal factor $e^{2\lambda}$. 
For the standard Euclidean coordinates $x_1,x_2$ on $D$ introduce the complex coordinates $z = x_1 + i x_2$ and  $\bar{z}$. Further, complexify the tangent space of $M$ and extend all tensors on it $\C$-linearly. For the remainder of this section we fix the following notation:
\begin{align*}
e_i &:= e^{-\lambda} \partial_{x_i} \phi & \\
\partial_z &:= \frac{\partial_{x_1} - i \partial_{x_2}}{2} &
\partial_{\bar{z}} & = \frac{\partial_{x_1} + i \partial_{x_2}}{2} \\
e_z &:= e^{-\lambda} \partial_{z} \phi = \frac{e_1 -i e_2}{2} &
e_{\bar{z}} &:= e^{-\lambda} \partial_{\bar{z}} \phi = \frac{e_1  + i e_2}{2}.
\end{align*}
By $\nabla^M_i$, $\nabla^M_z$ and $ \nabla^M_{\bar{z}}$ we mean $\nabla^M_{\partial_{x_i} \phi}$, $\nabla^M_{\partial_z \phi}$ and $\nabla^M_{\partial_{\bar{z}} \phi}$ respectively, and the following quantity can be seen as a complex version of the trace free second fundamental form.
\[ H_0 := 4 A(e_z,e_z) = A_{11} - A_{22} -2i A_{12} \]

\begin{lemma}[cf. {\cite[Lemma 3.3]{Mon_Riv_13_1}}]
\label{lemma_calc_multi_lin}
The pair $\{ e_1,e_2 \} $ is an orthonormal frame on $\Sigma$ and we can choose the orientation of $\{ e_1,e_2 \}$ such that $\{e_1,e_2,\nu \}$ with \[\nu = * e_1 \wedge e_2, \] is positively oriented. Here $*$, and $\wedge$ are the hogde star operator of $(M,g)$ and the wedge product respectively. Moreover, the following identities hold.
\begin{align*}
\langle e_z, e_z \rangle &= 0 = \langle e_{\bar{z}}, e_{\bar{z}} \rangle \\
 \langle e_z, e_{\bar{z}} \rangle &= \frac{1}{2} \\
\nabla^M_{z} \nu & = \frac{1}{2} H \p_z \phi + \frac{1}{2} H_0 \p_{\bar{z}} \phi \\
\end{align*}
\end{lemma}

\begin{theorem}[cf. {\cite[Theorem 3.1]{Mon_Riv_13_1}}]
\label{thm_div_form_general}
In the setting of this section the following identity holds:
\begin{align*}
4 e^{-2 \lambda} \Re \left(  \nabla^M_{\bar{z}} \left[  \partial_z H \nu + \frac{1}{2} H H_0 \partial_{\bar{z}} \phi \right] \right) = \Delta H \nu +H |\A|^2 \nu + 8 H \Re \left(  \tilde{\RR} e_z \right).
\end{align*}
where $\tilde{R} = g\left( \RR^M(e_{\bar{z}}, e_z)e_z,\nu \right)$. \\
In particular, if $\phi$ is of generalized Willmore type then we obtain the generalized Willmore equation in divergence form.
\begin{align}
\label{eq_gen_Will_div}
4  \Re \left(  \nabla^M_{\bar{z}} \left[  \partial_z H \nu + \frac{1}{2} H H_0 \partial_{\bar{z}} \phi \right] \right) = - e^{2 \lambda} F(\phi) \nu + 8 e^{2 \lambda} H \Re \left(  \tilde{\RR} e_z \right)
\end{align}
\end{theorem}

\noindent This theorem warrants the investigation of the vector field
\begin{align*} 
Y :=  H H_0 \partial_{\bar{z}} \phi + 2 \partial_z H \nu . 
\end{align*}
Note that it differs from $Y_f$, the one used in \cite[Equation 5.10]{Mon_Riv_13_1}, by setting $f=0 $ and multiplying by $-i$.

\noindent We introduce the space $(L^1 + W^{-1,2})(D)$ as the set of functions $f = f_1 + f_2$ with $f_1\in L^1(D)$ and $f_2 \in W^{-1,2}(D)$. 
It is equipped with the norm $\|f\|_{L^1+W^{-1,2}} = \inf_{f=f_1+f_2} \{ \|f_1\|_{L^1} +\|f_2\|_{W^{-1,2}} \}$.
By abuse of notation we will write $v \in X$ instead of $X(D, \R^n)$, also for (multi) vector fields $v$ and function spaces $X$.

Furthermore,  we assume that $\phi(D)$ is contained in a coordinate patch of $M$ such that we can trivialize $\phi^* TM \cong D \times \R^3$, where $\partial_1 \phi \mapsto b_1$, $\partial_2 \phi \mapsto b_2$, $\nu \mapsto b_3$ for $\{b_i\}$, the standard basis of $\R^3$. 
 In this trivialization $Y$  is a vector field in $(L^1 + W^{-1,2})(D, \C^3)$ since $A\in L^2(D)$. 

The next lemma establishes the existence of a potential $K$ of $Y$ and two more potentials $B_0$ and $B$ related to $K$. 
It corresponds to \cite[Lemma 6.1]{Mon_Riv_13_1}. 
Its proof is a direct application of the general constructions  \cite[Lemma A.1 and Lemma A.2]{Mon_Riv_13_1} and  depends on the fact that $\Re [\partial_{z} \phi \wedge (Y- 2 \nabla^M_z (H \nu))] = 0$ and $\Re[\langle \partial_{\bar{z}} \phi , Y \rangle ] =0 $. 
Further, it needs the regularity conditions $Y \in L^1+W^{-1,2}$ and $\Re(\nabla^M_{\bar{z}} Y) \in L^1 + W^{-1,2}$ both of which are true by the assumptions on $\phi$. 
In particular, it does not depend on the explicit shape of $\Re(\nabla^M_{\bar{z}} Y)$.

\begin{lemma}[cf. {\cite[Lemma 6.1]{Mon_Riv_13_1}}]
\label{lemma_potentials_construction}
Let $Y$ be the vector field from above. If $\phi$ is of generalized Willmore type, then
\begin{enumerate}
	\item[1)] there exists a complex vector field $K \in L^{q}$, with $\Im K \in W^{1,q}$, for every $q\in (1,2)$, that is the unique solution of
	\begin{align*}
\left\{ 
\begin{array}{ll}
\nabla^M_z K = i Y & \textup{in } D \\
 \Im K =0 & \textup{on } \partial D .
\end{array}
\right.
\end{align*}

\item[2)] There is a complex function $B_0 \in W^{1,q}$, with $\Im B_0 \in W^{2,q}$ for every $q\in (1,2)$ that solves
 	\begin{align*}
\left\{ 
\begin{array}{ll}
\partial_z B_0 = \langle \p_z \phi, \overline{K}  \rangle & \textup{in } D \\
 \Im B_0 =0 & \textup{on } \partial D .
\end{array}
\right.
\end{align*}

\item[3)]
There is a complex 2-vector field  $B \in W^{1,q}$, with $\Im B \in W^{2,q}$ for every $q\in (1,2)$ that solves
 	\begin{align*}
\left\{ 
\begin{array}{ll}
\nabla^M_z B = \p_z \phi \wedge \overline{K} + 2 i H \p_z \phi \wedge \nu & \textup{in } D \\
 \Im B =0 & \textup{on } \partial D .
\end{array}
\right.
\end{align*}
\end{enumerate}
\end{lemma}

\begin{definition}[cf. {\cite[Equation 1.34]{Mon_Riv_13_1}}]
For $u, v, w \in \Gamma(TM)$ define a contraction $\bullet$ of a vector field with a two vector field linarly on pure two vectors fields as follows: 
\[ u \bullet \left(v \wedge w \right) := \langle u, v \rangle w - \langle u , v \rangle w . \]
\end{definition}

The next lemma follows from Lemma \ref{lemma_potentials_construction} by  direct computation.
\begin{lemma}[cf. {\cite[Proposition 6.1]{Mon_Riv_13_1}]}]
\label{lemma_potentials_wente}
Let $\phi$ be of generalized Willmore type, then we have
\begin{align*}
\Delta \Re B &= *  \left[ \left(\nabla^M_2 \nu \right) \bullet \nabla^M_1 \Re B - \left(\nabla^M_1 \nu \right) \bullet \nabla^M_2 \Re B \right]  \\
& \hspace{0.5cm}-  \left(  \left( \p_1 \Re B_0  \right) \nabla^M_2 e_1\wedge e_2 - \left( \p_2 \Re B_0 \right) \nabla^M_1 e_1\wedge e_2 \right) + I , \\
\Delta \Re B_0 &= \langle \nabla^M_1 \Re B , \nabla^M_2 e_1\wedge e_2 \rangle -\langle \nabla^M_2 \Re B, \nabla^M_1 e_1\wedge e_2 \rangle  + G .
\end{align*}
Here $I$ and $G$ are functions in $L^q$, $q\in (1,2)$ which depend on $B, B_0, \Delta \Im B, \Delta \Im B_0$ derivatives of the metric and the second fundamental form $A$. 
\end{lemma}
 In coordinates the  system is of the form 
\[\Delta U^j = \p_1 E^j_k \p_2 U^k - \p_2 E^j_k  \p_1 U^k + \tilde{I}^j \]
for $U=\left( \Re (B_{ij}), \Re (B_0) \right)$, $\p_i U^j \in L^{q}$, $\p_i E^j_k \in L^2$ and where the $\tilde{I}^j \in L^q$, $q \in (1,2)$ are comprised of $I$ and $G$ and terms involving Christoffel symbols, $A$, $U^j$ and $\p_i U^j$.

\begin{theorem}
\label{thm_gen_will_smooth}
Let $\phi \in \FF(D,M)$ with conformal factor $e^{2\lambda}$, $\lambda \in L^\infty(D)$.
If $\phi$ is of  generalized Willmore type then $\phi$ is smooth.
\end{theorem}

\begin{proof}
In the first step we prove that $ (B, B_0) \in W_{loc}^{1,p}$ for a $p>2$.
This is done completely analogous to the proof of  \cite[Theorem 6.1]{Mon_Riv_13_1}. See also \cite{Sharp_Topping_Riv_equ} for a comprehensive treatment of the kind of PDE system that $ (B, B_0)$ solve.

In the second step we proceed differently. In particular, we adapt the bootstrap procedure between $H$ and $\phi$ to account for the function $F$ in Equation \eqref{eq_gen_Will_div}. 
The defining equation for $B$ reads
\[ 2i H \partial_z \phi \wedge \nu  =  \nabla^M_z B - \p_z \phi \wedge \overline{K}.  \]
Projecting this to $\partial_{\bar{z}}  \phi \wedge \nu$ and taking the imaginary part, yields 
\[H = e^{2\lambda} \Im\left( \left\langle \nabla^M_z B, \partial_{\bar{z}}  \phi \wedge \nu  \right\rangle \right) + \frac{1}{2} \langle \Im (K) , \nu \rangle,\]
 and hence $H \in L^{p}_{loc}$, as $\Im(K)\in L^q $ for all $q\in[1,\infty)$.
Since $\phi$ is conformal we have $\Delta_E \phi = e^{2\lambda} \vec{H}$, where $\Delta_E$ is the Euclidean Laplace operator, and hence $\phi \in W^{2,p}_{loc}$. 
In the following we retreat to $D_{1/2}(0)$, in order to drop the $loc$ subscript.

The generalized Willmore equation in divergence form reads
\[ 4 \Re\left(  \nabla^M_{\bar{z}}\left[ \p_z H \nu + \frac{1}{2} H H_0 \p_{\bar{z}} \phi \right]  \right)  = 8 \Re\left(  e^{2\lambda} H \tilde{R}e_z\right) - e^{2\lambda} F(\phi) \nu.  \]
In terms of a local frame $\{ b_\alpha \}$ of $M$ with Christoffel symbols $\Gamma$, $\nabla^M_{\bar{z}} \nabla^M_z \vec{H}$ can be expressed as 
\[ \nabla^M_{\bar{z}} \nabla^M_z \vec{H}  = \frac{1}{4} \Delta_E \vec{H} + \p_z \vec{H} \star \p_{\bar{z}} \phi \star\Gamma \star b + \p_{\bar{z}}(\vec{H} \star \p_z \phi \star \Gamma)  \star b + \vec{H} \star \p_z\phi \star \p_{\bar{z}} \phi \star \Gamma \star \Gamma  \star b .  \]
Here we employed the $ \star $ notation, that is $F \star G$ denotes a sum of contractions of $F$ and $G$.
Combining the last two equations we get an elliptic equation for $\vec{H}$ whose right hand side we control.
\begin{align}
\label{eq_h_local}
 \Delta_E \vec{H} &= 4 \Re \left(  \frac{1}{2} 
 \nabla^M_{\bar{z}}( H^2 \p_z \phi)  + 2 e^{2\lambda} H \tilde{R}e_z -\p_z \vec{H}   \star  \p_{\bar{z}} \phi  \star \Gamma  \star b \right. \\ 
& \hspace{0.7cm}  - \p_{\bar{z}}(\vec{H} \star \p_z \phi  \star \Gamma)  \star b  - \vec{H} \star \p_z\phi \star \p_{\bar{z}} \phi \star \Gamma \star \Gamma  \star b  \Big)   - e^{2\lambda} F(\phi) \nu \nonumber
\end{align}
By definition \ref{def_gen_will}, there is an $\epsilon >0 $ such that
	\begin{align*}
	e^{2\lambda} F(\phi) \nu &\in W^{k-1,l}, \ l= \frac{2p}{2+p} + \epsilon \ \ \textup{if } \phi \in W^{k+2,p} \cap W^{1,\infty} , \ p>2 , \ k\geq 0 .
	\end{align*}
Now suppose $\phi \in W^{k+2,p}$ for some $k\geq 0$, $p>2$ then the right hand side of \eqref{eq_h_local} is in $W^{k-1,l'}$, where  $l' = \min(l,p/2)$, if $k=0$ and $l' = \min(l, p)$ if $k>0$. Hence $\vec{H}\in W^{k+1,l'}$ and by  the equation $\Delta_E \phi = e^{2\lambda} \vec{H}$ we arrive at $\phi \in W^{k+3,l'}$. The following iteration implies the smoothness of $\phi$.

Let $p_0 := 2 + \delta $, for some  $0 < \delta < \epsilon/2 $ small enough such that $H\in L^{p_0}$, $\phi \in W^{2,p_0}$.
 Set $p_i :=  \frac{2l_{i-1}'}{2-l_{i-1}'}$ for $i\in \N$; $l_i := \frac{2 p_i}{ 2+ p_i} + \epsilon$ for $i \in \N_0$  and $l_i' := \min(l_i, p_i/2)$.
Since $ W^{1,l_i'} \hookrightarrow L^{p_{i+1}}$, we see that $\vec{H} \in W^{1,l_i'} $ for all $i\in \N$. 
As $p_i \to \infty$ and $l_i \to 2 + \epsilon $ we  eventually have that $\vec{H}\in W^{1,p_0}$.
Now we may iterate again for the higher derivatives.

\end{proof}

\bibliographystyle{hplain} 
\bibliography{thesis_bib} 

\end{document}